\documentclass[11pt]{amsart}
\usepackage[UKenglish]{babel}
\usepackage[utf8]{inputenc} % deutsche Umlaute
\usepackage{amsthm} % Satz/Beweis-Strukturen
\usepackage{amsfonts} % Die Mengen IN, Z, Q, IR, C, ...
\usepackage{mathtools}
\usepackage{enumitem} \setlist[enumerate]{label={\upshape(\arabic*)}}
\usepackage{amssymb}
\usepackage[numbers]{natbib}
\usepackage{url}
\usepackage{mathrsfs}
\usepackage{geometry}
\usepackage{bbm}
\usepackage{tikz-cd}
\usepackage{tikz}
\usepackage{color}
\usepackage{framed}
\usepackage{hyperref}
\hypersetup{
    colorlinks=true, % make the links colored
    linkcolor=blue, % color TOC links in blue
    urlcolor=red, % color URLs in red
    citecolor=[rgb]{0,0.7,0},
    linktoc=all % 'all' will create links for everything in the TOC
}

\theoremstyle{definition}
\newtheorem{defn}{Definition}[section]
\newtheorem{prop}[defn]{Proposition}

\newtheorem{thm}[defn]{Theorem}

\newtheorem{conj}[defn]{Conjecture}

\newcommand{\R}{\mathbb{R}}

\newcommand{\Z}{\mathbb{Z}}

\def\1{\mathbf{1}}

\newcommand{\N}{{\mathbb{N}}}

\newcommand{\hs}{h^*}

\title[Examples of IDP polytopes with non-log-concave $h^*$-vector]{Examples of IDP lattice polytopes with non-log-concave $h^*$-vector}
\author[J.~Hofscheier]{Johannes Hofscheier}
\author[V.~Kurylenko]{Vadym Kurylenko}
\author[B.~Nill]{Benjamin Nill}

\address[J.~Hofscheier]{School of Mathematical Sciences\\University of Nottingham\\ Nottingham\\NG7 2RD\\UK}
\email{johannes.hofscheier@nottingham.ac.uk}

\address[V.~Kurylenko, B.~Nill]{Faculty of Mathematics, Otto-von-Guericke-Universit\"at Magdeburg, Universit\"atsplatz 2, 39106 Magdeburg, Germany.}
\email{\{vadym.kurylenko, benjamin.nill\}@ovgu.de}


\begin{document}

\setlength{\parindent}{0pt}
\begin{abstract}
    Lattice polytopes are called {\em IDP polytopes} if they have the integer decomposition property, i.e., any lattice point in a $k$th dilation is a sum of $k$ lattice points in the polytope.
    It is a long-standing conjecture whether the numerator of the Ehrhart series of an IDP polytope, called the $h^*$-polynomial, has a unimodal coefficient vector.
    In this preliminary report on research in progress we present examples showing that $h^*$-vectors of IDP polytopes do not have to be log-concave. This answers a question of Luis Ferroni and Akihiro Higashitani. 
    
    As this is an ongoing project, this paper will be updated with more details and examples in the near future.
\end{abstract}

\maketitle

\thispagestyle{empty}

\section{Counterexamples to the log-concavity question}

In this note we report on examples found in an ongoing project to search for counterexamples to questions and conjectures in Ehrhart theory employing methods and heuristics from machine learning.

\medskip

Let us recall the main definitions. We refer to \cite{beckbook} for references. A {\em lattice polytope} $P \subset \R^n$ is the convex hull of finitely many lattice points, i.e., elements of $\Z^n$. The {\em Ehrhart polynomial} of $P$ is the unique polynomial $E_P(t)$ with $E_P(k) = |k P \cap \Z^n|$ for $k \in \N_{\ge 1}$. The {\em $\hs$-polynomial} $\hs_P(t)$ of $P$ is defined as the numerator of its {\em Ehrhart generating series}: 
\[\sum_{k=0}^\infty E_P(k) t^k = \frac{\hs_P(t)}{(1-t)^{d+1}},\]
where $d$ denotes the dimension of $P$. For $\hs_P(t)=\sum_{i=0}^d \hs_i t^i$, we call $(\hs_0, \ldots, \hs_d)$ the {\em $\hs$-vector} of $P$. 

\smallskip

The following properties of lattice polytopes have been intensively studied over the previous years. We refer to the excellent survey papers \cite{Braun-Survey} and \cite{Ferroni}. An {\em IDP polytope} (an abbreviation for lattice polytopes with the \emph{Integer Decomposition Property}) is a lattice polytope $P$ such that for every $k \in \N_{\ge 2}$ each lattice point in $kP$ is a sum of $k$ lattice points in $P$. We say, the $\hs$-vector of $P$ is {\em unimodal} if $\hs_0 \le \hs_1 \le \cdots \le \hs_j\ge \hs_{j+1} \ge \cdots \ge \hs_d$ for some $j \in \{0, \ldots, d\}$. It is called {\em log-concave} if $\hs_{i-1} \hs_{i+1} \le (\hs_i)^2$ for every $i \in \{1, \ldots, d-1\}$. 

\smallskip

The following conjecture is currently considered as the main open question in Ehrhart theory. It is explicitly stated in \cite[Conjecture~1.1]{Ferroni} and \cite[Question~1.1]{Van}.

\begin{conj} Every IDP polytope has a unimodal $\hs$-vector.\label{unimodality}  
\end{conj}

As $\hs$-vectors of IDP lattice polytopes do not have internal zeroes, log-concavity of the $\hs$-vector implies its unimodality. It is natural to wonder whether one can replace unimodal by log-concave in the previous conjecture. This is Question~3.9(b) in the paper by Ferroni and Higashitani \cite{Ferroni}. In the following we show that this does not hold. We remark that in the more general commutative-algebraic setting of $h$-vectors of standard Cohen-Macaulay domains, log-concavity had already been disproved in 1993 \cite{Niesi}.

\begin{thm}
Let $P \subset \R^7$ be the lattice polytope given by the convex hull of $e_1, \ldots, e_7$ and the following two vectors
\[(1,-1,-1,-1,-1,0,0)\] and 
\[(-1,-1,0,0,0,-1,-1).\]
Then $P$ has the following properties:
\begin{enumerate}
    \item $P$ is a $7$-dimensional IDP polytope with $9$ vertices.
    \item All coefficients of the Ehrhart polynomial $E_P(t)$ of $P$ are positive.
    \item The $\hs$-vector of $P$ is not log-concave:
    \[(\hs_0, \ldots, \hs_7) = (1,2,3,4,5,3,2,1)\]
\end{enumerate}
\end{thm}

The reader may notice that the vertices of the previous example have coefficients in $\{1,0,-1\}$. Our next example shows that in higher dimensions one can even find examples that are among the very special class of \emph{0/1-polytopes}. 

\begin{thm}
Let $P \subset \R^{12} $ be the lattice polytope given by the convex hull of $0, e_2, \ldots, e_{12}$ and the following three vectors in $\R^{12}$
\[(1, 0, 0, 0, 1, 1, 1, 1, 0, 0, 0, 0),\] 
\[(1, 1, 1, 1, 0, 0, 0, 0, 0, 0, 0, 0),\]
\[(1, 0, 0, 0, 0, 0, 0, 0, 1, 1, 1, 1).\]
Then $P$ has the following properties:
\begin{enumerate}
\item $P$ is a $12$-dimensional IDP 0/1-polytope with $15$ vertices.
\item Every triangulation of $P$ is unimodular.
\item $P$ has no quadratic triangulation (i.e., regular, unimodular, flag triangulation)
\item All coefficients of the Ehrhart polynomial $E_P(t)$ of $P$ are positive.
\item The $\hs$-vector of $P$ is not log-concave:
\[(\hs_0, \ldots, \hs_{12}) = (1,2,3,4,5,3,2,1,0,0,0,0,0)\]
\end{enumerate}
\label{12-bsp}
\end{thm}

\medskip
As the reader notices from (2), the previous example belongs to an extremely special kind of lattice polytopes: {\em unimodular polytopes} (see e.g. \cite{Unimodular}) in the very top of the hierarchy of lattice polytopes (see e.g. \cite{triangulations}). The reason the project members were looking for counterexamples in this class is due to the belief of the third author that if Conjecture~\ref{unimodality} is wrong it should fail even for very special classes. 
        
The reader might also be surprised by the fact that the above example has the small number of $d+3$ vertices.
This has been the result of applying machine learning techniques, namely, reinforcement learning, to decrease a measure of log-concavity by locally modifying polytopes. In all existing runs this led to IDP polytopes with $d+k$ vertices, where $k$ was a small number, which led us to look at polytopes with this number of vertices more closely.

\medskip

Finally, let us note that the lattice polytope in Theorem~\ref{12-bsp} comes from a directed graph. Let us recall (see e.g., \cite{triangulations}) that given a directed graph $D = (V,A)$ one can associate its {\em root polytope} or {\em arc polytope} as the convex hull of the columns of its incidence matrix (which is a $|V| \times |A|$-matrix with entries in $\{1,0,-1\}$). 

\begin{figure}[!ht]
    \centering
    \begin{tikzpicture}[
        every node/.style={ draw, circle, minimum size=6mm },
        ->, >=Stealth, thick
    ]
    
    \node (1) at (-1.75,-3) {1};
    \node (2) at (-1.75,-2) {2};
    \node (3) at (-1.75,-1) {3};
    \node (4) at (-1.75,0.00) {4};
    \node (5) at (-1.75,1) {5};
    \node (6) at (-1.75,2) {6};
    \node (7) at (-1.75,3) {7};
    \node (8) at (1.75,-3) {8};
    \node (9) at (1.75,-2) {9};
    \node (10) at (1.75,-1) {10};
    \node (11) at (1.75,0.00) {11};
    \node (12) at (1.75,1) {12};
    \node (13) at (1.75,2) {13};
    \node (14) at (1.75,3) {14};
    
    \draw (1) to (8);
    \draw (1) to (12);
    \draw (1) to (14);
    \draw (2) to (8);
    \draw (3) to (9);
    \draw (4) to (11);
    \draw (5) to (12);
    \draw (6) to (13);
    \draw (7) to (14);
    \draw (9) to (2);
    \draw (10) to (3);
    \draw (10) to (4);
    \draw (10) to (6);
    \draw (11) to (5);
    \draw (13) to (7);
    
    \end{tikzpicture}
    \caption{A bipartite graph realising the \(12\)-dimensional example.}
    \label{graph}
\end{figure}
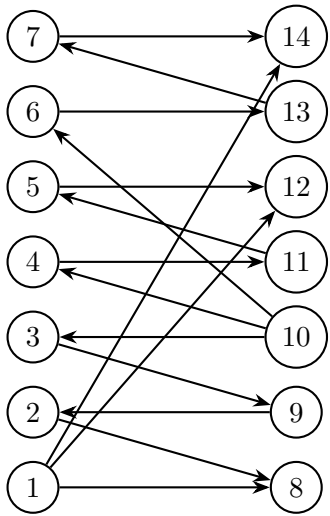

\begin{prop}
    The lattice polytope from Theorem~\ref{12-bsp} is unimodularly equivalent to the arc polytope of the directed bipartite graph in Figure~\ref{graph}.
\end{prop}

All properties in these results can be verified using software packages like Polymake, Normaliz, Magma, or SageMath.        
\medskip

We will add more details on the methods used to find these examples in a later version of this paper. It will also include more examples and observations, discussions on other interesting properties, and a more systematic report on our ongoing search for counterexamples.

\subsection*{Acknowledgment}
This work began during the first author's visit to Otto von Guericke University Magdeburg, whose hospitality is gratefully acknowledged.
The second and third authors are funded by the Deutsche Forschungsgemeinschaft (DFG, German Research Foundation) – 539867500 as part of the research priority program Combinatorial Synergies.
We thank Luis Ferroni for his feedback and crosschecking the example in Theorem~\ref{12-bsp}.
We are grateful to Christian Haase for his interest in this project.
\bibliographystyle{plain}
\bibliography{ref}

\end{document}